\documentclass[10pt]{amsart}
\input{epsf}
\usepackage{amssymb,latexsym}
\topskip  \numberwithin{equation}{section}
\newtheorem{theorem}{Theorem}[section]
\newtheorem{proposition}[theorem]{Proposition}
\newtheorem{corollary}[theorem]{Corollary}
\newtheorem{definition}[theorem]{Definition}

\newtheorem{example}[theorem]{Example}

\begin{document}

\pagenumbering{arabic}
\pagestyle{headings}
\def\sof{\hfill\rule{2mm}{2mm}}
\def\llim{\lim_{n\rightarrow\infty}}

\title{Partially Ordered generalized patterns and $k$-ary words}
\author{Sergey Kitaev and Toufik Mansour}
\maketitle
\begin{center}
Matematik, Chalmers tekniska h\"ogskola och G\"oteborgs universitet,\\
S-412~96 G\"oteborg, Sweden

{\tt kitaev@math.chalmers.se, toufik@math.chalmers.se} 
\end{center}
\def\P{POGP}
\def\A{\mathcal{A}}
\def\SS{\frak S}
\def\Ps{POGPs}
\def\mn{\mbox{-}}
\def\newop#1{\expandafter\def\csname #1\endcsname{\mathop{\rm #1}\nolimits}}
\newop{MND}
\section*{Abstract}
Recently, Kitaev~\cite{Ki2} introduced partially ordered generalized patterns (POGPs) in the symmetric group, which further generalize the generalized
permutation patterns introduced by Babson and Steingr{\'\i}msson~\cite{BS}. A POGP $p$ is a GP some of whose letters are incomparable. In
this paper, we study the generating functions (g.f.) for the number of
$k$-ary words avoiding some POGPs. We give analogues, extend and generalize several known results, as well as get some new results. In particular, we give the g.f. for the entire distribution of the maximum number of non-overlapping occurrences of a pattern $p$ with no hyphens (that allowed to have repetition of letters), provided we know the g.f. for the number of $k$-ary words that avoid $p$. 

\thispagestyle{empty} 
\section{Introduction}
Let $[k]^n$ denote the set of all the words of length $n$ over the (totally ordered) alphabet $[k]=\{ 1,2,\dots,k \}$. We call these words by $n$-long
$k$-ary words. A \emph{generalized pattern} $\tau$ is a word in $[\ell]^m$ (possibly with hyphens between some letters) that contains each letter from $[\ell]$ (possibly with repetitions). We say that the word $\sigma\in[k]^n$
\emph{contains} a generalized pattern $\tau$, if $\sigma$ contains
a subsequence isomorphic to $\tau$ in which the entries
corresponding to consecutive entries of $\tau$, which are not separated by a
hyphen, must be adjacent. Otherwise, we say that $\sigma$
\emph{avoids} $\tau$ and write $\sigma\in [k]^n(\tau)$. Thus,
$[k]^n(\tau)$ denotes the set of all the words in $[k]^n$ that avoid $\tau$. Moreover, if $P$ is a set of generalized patterns then $|[k]^n(P)|$ denotes the set all the words in $[k]^n$ that avoid each pattern from $P$ simultaneously. 

\begin{example}
A word $\pi=a_1a_2\dots a_n$ avoids the pattern $13\mn2$ if $\pi$ has no
subsequence $a_ia_{i+1}a_j$ with $j>i+1$ and $a_i<a_j<a_{i+1}$. Also, $\pi$ avoids the pattern $121$ if it has no subword $a_ia_{i+1}a_{i+2}$ such that $a_i=a_{i+2}<a_{i+1}$. 
\end{example}

{\em Classical patterns} are generalized patterns with all possible hyphens (say, $2\mn1\mn3\mn4$), in other words, those that place no
adjacency requirements on $\sigma$. The first case of classical
patterns studied was that of permutations avoiding a pattern of length 3 in $\SS_3$. Knuth \cite{Knuth} found that, for any
$\tau\in\SS_3$, $|\SS_n(\tau)|=C_n$, the $n$th Catalan number.
Later, Simion and Schmidt \cite{SS} determined the number
$|\SS_n(P)|$ of permutations in $\SS_n$ simultaneously avoiding
any given set of patterns $P\subseteq\SS_3$. Burstein
\cite{Burstein} extended this to $|[k]^n(P)|$ with
$P\subseteq\SS_3$. Burstein and Mansour \cite{BM1} considered
forbidden patterns with repeated letters. Also, Burstein and
Mansour \cite{BM2,BM3} considered forbidden generalized patterns
with repeated letters.

{\em Generalized permutation patterns} were introduced by Babson and
Steingr\'{\i}msson \cite{BS} with the purpose of the study of
Mahonian statistics. Claesson \cite{Claesson} and Claesson
and Mansour \cite{CM} considered the number of permutations
avoiding one or two generalized patterns with one hyphen. Kitaev 
\cite{Ki1} examined the number of $|\SS_n(P)|$ of
permutations in $\SS_n$ simultaneously avoiding any set of
generalized patterns with no hyphens. Besides, Kitaev~\cite{Ki2}
introduced a further generalization of the generalized permutation patterns namely
{\em partially ordered generalized patterns}.

In this paper we introduce a further generalization of the generalized
patterns namely {\em partially ordered generalized patterns in words {\rm(}\Ps{\rm)}}, 
which is an analogue of \Ps\, in permutations \cite{Ki2}. 
A \P\, is a generalized pattern some of whose letters are incomparable. 
For example, if we
write $\tau=1\mn1'2'$, then we mean that in occurrence of $\tau$
in a word $\sigma\in[k]^n$ the letter corresponding to the $1$ in
$\tau$ can be either larger, smaller, or equal to the letters
corresponding to $1'2'$. Thus, the word $113425\in[5]^6$ contains
seven occurrence of $\tau$, namely $113$, $134$ twice, $125$ twice, $325$, 
and $425$.

Following~\cite{Ki2}, we consider two particular classes of \Ps\, -- 
{\em shuffle patterns} and {\em multi-patterns}, which allows us to
give an analogue for all the main results of~\cite{Ki2} for $k$-ary
words. A multi-pattern is of the form $\tau=\tau^0\mn\tau^1\mn\cdots\mn\tau^s$ 
and a shuffle pattern of the form $\tau=\tau^0\mn a_1\mn\tau^1\mn a_2\mn\cdots\mn\tau^{s-1}\mn a_s\mn\tau^s$, where for any $i$ and $j$, the letter $a_i$ is greater than any letter of $\tau^j$ and for any $i\neq j$ each letter of $\tau^i$ is incomparable with any letter of $\tau^j$. These patterns are investigated in Sections~\ref{sec3} and~\ref{sec4}.

Let $\tau=\tau^0\mn\tau^1\mn\cdots\mn\tau^s$ be an arbitrary multi-pattern and let $A_{\tau^i}(x;k)$ be the generating function (g.f.) for the number of words in $k$-letter alphabet that avoid $\tau^i$ for each $i$. In Theorem~\ref{coaa} we find the g.f., in terms of the $A_{\tau^i}(x;k)$, for the number of $k$-ary words that avoid $\tau$. In particular, this allows us to find the g.f. for the entire {\em distribution} of the maximum number of non-overlapping occurrences of a pattern $\tau$ with no hyphens, if we only know the g.f. for the number of $k$-ary words that avoid $\tau$. Thus, in order to apply our results in what follows we need to know how many $k$-ary words avoid a given ordinary generalized pattern with no hyphens. This question was examined, for instance, in~\cite[Sections 2 and 3]{BM1},~\cite[Section 3]{BM2} and~\cite[Section 3.3]{BM3}. 

\section{Definitions and Preliminaries}
A {\em partially ordered generalized pattern {\rm(}\P{\rm)}} is a
generalized pattern where some of the letters can be incomparable.

\begin{example}
The simplest non-trivial example of a \P\, that differs from the
ordinary generalized patterns is $\tau=1'\mn 2\mn 1''$, where the
second letters is the greatest one and the first and the last
letters are incomparable to each other. The word $\sigma=31421$
has five occurrences of $\tau$, namely $342$, $341$, $142$, $141$,
and $121$.
\end{example}

Let $A_{\tau}(x;k)=\sum_{n\geq0}a_\tau(n;k)x^n$ denote the {\em generating function {\rm(}g.f.{\rm)}} for the numbers $a_\tau(n;k)$ of words in $[k]^n$ avoiding the pattern~$\tau$. For $\tau=1'\mn2\mn1''$, we have 
\begin{equation}
A_{1'\mn2\mn1''}(x;k) = \frac{1}{(1-x)^{2k-1}}-\sum_{j=1}^{k-1}
\frac{x}{(1-x)^{2j}}.\label{eqex1}
\end{equation}
Indeed, if $\sigma\in [k]^n$ avoids $\tau$, and $\sigma$ contains $s>0$ copies of
the letter $k$, then the letters $k$ appear as leftmost or rightmost letters of
$\sigma$. If $\sigma$ contains no $k$ then $\sigma\in[k-1]^n$. So, for all $n\geq0$, we have
$$a_\tau(n;k)=a_\tau(n;k-1)+2a_\tau(n-1;k-1)+3a_\tau(n-2;k-1)+\cdots+(n+1)a_\tau(0;k-1),$$
since there are $(i+1)a_\tau(n-i;k-1)$ possibilities to place $i$ letters $k$ into $\sigma$, for $0\leq i \leq n$. Hence, for all $n\geq2$,
$$a_\tau(n;k)-2a_\tau(n-1;k)+a_\tau(n-2;k)=a_\tau(n;k-1),$$
together with $a_\tau(0,k)=1$ and $a_\tau(1,k)=k$. Multiplying
both sides of the recurrence above with $x^n$ and summing over all
$n\geq2$, we get Equation~\ref{eqex1}.

\begin{definition}
If the number of words in $[k]^n$, for each $n$, that avoid a \P\, $\tau$ 
is equal to the number of words that avoid a \P\, $\phi$, then $\tau$ and 
$\phi$ are said to be {\em equivalent} and we write $\tau\equiv\phi$.
\end{definition}

The {\em reverse $R(\sigma)$} of a word
$\sigma=\sigma_1\sigma_2\ldots\sigma_n$ is the word
$\sigma_n\ldots\sigma_2\sigma_1$. The {\em complement $C(\sigma)$}
is the word $\theta=\theta_1\theta_2\ldots\theta_n$ where
$\theta_i=k+1-\sigma_i$ for all $i=1,2,\ldots,n$. For example, if
$\sigma=123331\in[3]^6$, then $R(\sigma)=133321$,
$C(\sigma)=321113$, and $R(C(\sigma))=311123$. We call these bijections 
of $[k]^n$ to itself {\em trivial}. For example, the number of words that avoid the
pattern $12\mn2$ is the same as the number of words that avoid the
patterns $2\mn21$, $1\mn12$, and $21\mn1$, respectively.

Following~\cite{Ki2}, it is convenient to introduce the following definition.

\begin{definition}
Let $\tau$ be a generalized pattern without hyphens. A word
$\sigma$ {\em quasi-avoids} $\tau$ if $\sigma$ has exactly one
occurrence of $\tau$ and this occurrence consists of the $|\tau|$ rightmost
letters of $\sigma$, where $|\tau|$ denotes the number of letters in $\tau$.
\end{definition}

For example, the word $5112234$ quasi-avoids the pattern $1123$,
whereas the words $5223411$ and $1123345$ do not. 

\begin{proposition}\label{praa}
Let $\tau$ be a non-empty generalized pattern with no hyphens. Let
$A_\tau^*(x;k)$ denote the g.f. for the number of words in $[k]^n$ that 
quasi-avoid $\tau$. Then
\begin{equation}
A_\tau^*(x;k)=(kx-1)A_\tau(x;k)+1.
\end{equation}
\end{proposition}
\begin{proof}
Using the similar arguments as those in the proof of \cite[Proposition~4]{Ki2}, 
we get that, for $n\geq 1$,
$$a_\tau^*(n;k)=ka_\tau(n-1;k)-a_\tau(n;k),$$
where $a_\tau^*(n;k)$ denotes the number of words in $[k]^n$ that quasi-avoid 
$\tau$. Multiplying both sides of the last equality by $x^n$ and summing over 
all natural numbers $n$, we get the desired result.
\end{proof}

\begin{definition}\label{defmult}
Suppose $\{\tau^0,\tau^1,\dots,\tau^s\}$ is a set of generalized
patterns with no hyphens and
$$\tau=\tau^0\mn\tau^1\mn\cdots\mn\tau^s,$$ where each letter of
$\tau^i$ is incomparable with any letter of $\tau^j$ whenever
$i\neq j$. We call such \Ps\, {\em multi-patterns}.
\end{definition}

\begin{definition}\label{defshuf}
Suppose $\{\tau^0,\tau^1,\dots,\tau^s\}$ is a set of generalized
patterns with no hyphens and $a_1a_2\ldots a_s$ is a word of $s$
letters. We define a {\em shuffle} pattern to be a pattern  of the
form
$$\tau=\tau^0\mn a_1\mn\tau^1\mn a_2\mn\cdots\mn\tau^{s-1}\mn a_s\mn\tau^s,$$
where each letter of $\tau^i$ is incomparable with any letter of
$\tau^j$ whenever $i\neq j$, and the letter $a_i$ is greater than
any letter of $\tau^j$ for any $i$ and $j$.
\end{definition}

For example, $1'\mn2\mn1''$ is a shuffle pattern, and $1'\mn1''$
is a multi-patterns. From definitions, we obtain that we can get a
multi-pattern from a shuffle pattern by removing all the letters
$a_i$. 

There is a connection between multi-avoidance of the generalized patterns 
and the \Ps. In particular, to avoid $1'\mn2\mn1''$ is the same as to
avoid simultaneously the patterns $1\mn2\mn1$, $1\mn3\mn2$, and
$2\mn3\mn1$. A straightforward argument leads to the following
proposition.

\begin{proposition}
For any \P\, $\tau$ there exists a set $T$ of generalized patterns
such that a word $\sigma$ avoids $\tau$ if and only if $\sigma$
avoids all the patterns in $T$.
\end{proposition}

For example, if $\tau=1'2'\mn3\mn1''$, then to avoid $\tau$ is the
same to avoid $5$ patterns, $12\mn3\mn1$, $12\mn3\mn2$,
$12\mn4\mn3$, $13\mn4\mn2$, and $23\mn4\mn1$. Moreover, the following proposition 
holds:

\begin{proposition} Suppose $\tau=\tau_1\mn a \mn\tau_2$ {\rm(}resp. $\phi = \phi_1\mn\phi_2${\rm)} is a 
shuffle pattern 
{\rm(}resp. a multi-pattern{\rm)} such that $\tau_1,\ \phi_1 \in [r_1]^{\ell_1}$, $\tau_2,\ \phi_2 \in [r_2]^{\ell_2}$ 
and each letter of $[r_1]$ is incomparable with any letter of $[r_2]$. 
Also, without lose the generality, suppose $r_1\geq r_2$. Then to avoid $\tau$ {\rm(}resp. $\phi${\rm)} is 
the same as to avoid $\displaystyle\sum_{i=0}^{r_2}{r_1 \choose i}{r_2 \choose i}{r_1+r_2-i \choose r_1}$ 
generalazed patterns. In particular, the number of generalized patterns does not depend on the lengths $\ell_1$ and $\ell_2$. 
\end{proposition}
\begin{proof}
Obviously, to prove the statement, we need to find the number of ways to make a total order on $[r_1] \cup [r_2]$ (the letter $a$ does not play any roll, since it is always the greatest letter). Any total order on $[r_1] \cup [r_2]$ is an alphabet that can consist of $r_1+r_2-i$ letters, where $i$ is the number of letters in $[r_2]$ that supposed to coincide with some letters in $[r_1]$. Clearly, $0\leq i \leq r_2$ and we can choose coinciding letters in ${r_1 \choose i}{r_2 \choose i}$ ways. Now, after choosing the coinciding letters, we can make a total order in ${r_1+r_2-i \choose r_1}$ ways, which is given by~\cite[Theorem 8]{Ki2}.      
\end{proof}
\section{The shuffle pattern}\label{sec3}
We recall that according to Definition~\ref{defshuf}, a shuffle pattern is a pattern of the form
$$\tau=\tau^0\mn a_1\mn\tau^1\mn a_2\mn\cdots\mn\tau^{s-1}\mn a_s\mn\tau^s,$$ where 
$\{\tau^0,\tau^1,\dots,\tau^s\}$ is a set of generalized patterns with no hyphens, 
$a_1a_2\ldots a_s$ is a word of $s$ letters, for any $i$ and $j$ the letter $a_i$ is greater than any letter of 
$\tau^j$ and for any $i \neq j$ each letter of $\tau^i$ is incomparable with any letter of $\tau^j$.
 
Let us consider the shuffle pattern $\phi = \tau\mn \ell \mn\tau$, where $\ell$ is the greatest letter in $\phi$ and letters each letter in the left $\tau$ is incomparable with any letter in the right $\tau$.

\begin{theorem}\label{thsa}
Let $\phi$ be the shuffle pattern $\tau\mn \ell \mn\tau$ described above. Then for all $k\geq \ell$,
$$A_\phi(x;k)=\frac{1}{(1-xA_\tau(x;k-1))^2}\biggl(A_\phi(x;k-1)-xA_\tau^2(x;k-1)\biggr).$$
\end{theorem}
\begin{proof}
We show how to get a recurrence relation on $k$ for $A_\phi(x;k)$, which is the g.f. for the number of words in $[k]^n(\phi)$. 
Suppose $\sigma\in [k]^n(\phi)$ is such that it contains exactly $d$ 
copies of the letter $k$. If $d=0$ then the g.f. for the number of such words is $A_\phi(x;k-1)$. Assume that $d\geq 1$. Clearly, 
$\sigma$ can be written in the following form:
$$\sigma=\sigma^0 k \sigma^1 k \cdots k \sigma^{d},$$
where $\sigma^j$ is a $\phi$-avoiding word on $k-1$ letters, for $j=0,1, \ldots, d$. There 
are two possibilities: either $\sigma^j$ avoids $\tau$ for all $j$, or there exists $j_0$ 
such that $\sigma^{j_0}$ contains $\tau$ and for any $j\neq j_0$, the word $\sigma^j$ avoids~$\tau$. In the first case, the number of such words is given by the g.f. $x^dA_\tau^{d+1}(x;k-1)$, whereas in the second case, by $(d+1)x^dA_\tau^d(x;k-1)(A_\phi(x;k-1)-A_\tau(x;k-1))$. In the last expression, the multiple $(d+1)$ is the number of ways to choose $j$, such that $\sigma^j$ has an occurrence of $\tau$, and $A_\phi(x;k-1)-A_\tau(x;k-1)$ is the g.f. for the number of words avoiding $\phi$ and containing $\tau$.   

Therefore,
$$A_\phi(x;k)=A_\phi(x;k-1)+\sum_{d\geq 1}(d+1)x^dA_\tau^d(x;k-1)A_\phi(x;k-1)-
\sum_{d\geq 1} dx^dA_\tau^{d+1}(x;k-1),$$ equivalently,
$$A_\phi(x;k)=A_\phi(x;k-1)+A_\phi(x;k-1)\frac{2xA_\tau(x;k-1)-x^2A_\tau^2(x;k-1)}{(1-xA_\tau(x;k-1))^2}
-\frac{xA_\tau^2(x;k-1)}{(1-xA_\tau(x;k-1))^2}.$$ The rest is easy
to check.
\end{proof}

\begin{example}
Let $\phi=1'\mn2\mn1''$. Here $\tau=1$, so $A_\tau(x;k)=1$ for all $k\geq 1$, since 
only the empty word avoids $\tau$. Hence, according to Theorem~\ref{thsa}, we have
$$A_\phi(x;k)=\frac{A_\phi(x;k-1)-x}{(1-x)^2},$$
which together with $A_\phi(x;1)=\frac{1}{1-x}$ {\rm(}for any $n$ only the word $\underbrace{11\ldots 1}_{n\ times}$ avoids $\phi${\rm)} gives Equation~\ref{eqex1}.
\end{example}

More generally, we consider a shuffle pattern of the form $\tau^0\mn \ell\mn \tau^1$, where $\ell$ 
is the greatest element of the pattern.

\begin{theorem}\label{thsb}
Let $\phi$ be the shuffle pattern $\tau\mn \ell \mn\nu$. Then for all
$k\geq \ell$,
$$A_\phi(x;k)=\frac{1}{(1-xA_\tau(x;k-1))(1-xA_\nu(x;k-1))}\biggl(A_\phi(x;k-1)-xA_\tau(x;k-1)A_\nu(x;k-1)\biggr).$$
\end{theorem}
\begin{proof}
We proceed as in the proof of Theorem~\ref{thsa}. Suppose $\sigma\in [k]^n(\phi)$ is such that it contains exactly $d$ 
copies of the letter $k$. If $d=0$ then the g.f. for the number of such words is $A_\phi(x;k-1)$. Assume that $d\geq 1$. Clearly, 
$\sigma$ can be written in the following form:
$$\sigma=\sigma^0 k \sigma^1 k \cdots k \sigma^{d},$$
where $\sigma^j$ is a $\phi$-avoiding word on $k-1$ letters, for $j=0,1, \ldots, d$. There 
are two possibilities: either $\sigma^j$ avoids $\tau$ for all $j$, or there exists $j_0$ 
such that $\sigma^{j_0}$ contains $\tau$, $\sigma^j$ avoids $\tau$ for all $j=0,1, \ldots, j_0-1$ and $\sigma^j$ avoids~$\nu$ for any $j=j_0+1,\ldots,d$. In the first case, the number of such words is given by the g.f. $x^dA_\tau^{d+1}(x;k-1)$. In the second case, we have
$$x^d\sum_{j=0}^dA_\tau^j(x;k-1)A_\nu^{d-j}(x;k-1)(A_\phi(x;k-1)-A_\tau(x;k-1)).$$
Therefore, we get
$$\begin{array}{l}
A_\phi(x;k)=A_\phi(x;k-1)+A_\phi(x;k-1)\sum\limits_{d\geq1}
x^d\sum\limits_{j=0}^dA_\tau^j(x;k-1)A_\nu^{d-j}(x;k-1)\\
\qquad\qquad\qquad\qquad\qquad\qquad\qquad\qquad-\sum\limits_{d\geq1}
x^d\sum\limits_{j=1}^d
A_\tau^j(x;k-1)A_\nu^{d+1-j}(x;k-1),\end{array}$$ equivalently,
$$A_\phi(x;k)=(A_\phi(x;k-1)-xA_\tau(x;k-1)A_\nu(x;k-1))\sum_{d\geq0}
x^d\sum_{j=0}^dA_\tau^j(x;k-1)A_\nu^{d-j}(x;k-1).$$ Hence, using
the identity $\displaystyle\sum_{n\geq0} x^n\sum_{j=0}^n
p^jq^{n-j}=\frac{1}{(1-xp)(1-xq)}$ we get the desired result.
\end{proof}

We now give two corollaries to Theorem~\ref{thsb}.

\begin{corollary}\label{thsc}
Let $\phi=\tau^0\mn \ell \mn\tau^1$ be a shuffle pattern, and let
$f(\phi)=f_1(\tau^0)\mn \ell \mn f_2(\tau^1)$, where $f_1$ and $f_2$
are any trivial bijections. Then $\phi\equiv f(\phi)$.
\end{corollary}
\begin{proof}
Using Theorem~\ref{thsb}, and the fact that the number of words in
$[k]^n$ avoiding $\tau$ (resp. $\nu$) and $f_1(\tau)$
(resp. $f_2(\nu)$) have the same generating functions, we get the desired result.
\end{proof}

\begin{corollary}\label{thsd}
For any shuffle pattern $\tau\mn \ell \mn\nu$, we have
$$\tau\mn \ell \mn\nu\equiv \nu\mn \ell \mn\tau.$$
\end{corollary}
\begin{proof}
Corollary~\ref{thsc} yields that the shuffle pattern $\tau\mn \ell \mn\nu$ is
equivalent to the pattern $\tau\mn \ell \mn R(\nu)$, which is equivalent to
the pattern $R(\tau\mn \ell \mn R(\nu))=\nu\mn \ell \mn R(\tau)$. Finally, we use 
Corollary~\ref{thsc} one more time to get the desired result.
\end{proof}
\section{The multi-patterns}\label{sec4}
We recall that according to Definition~\ref{defmult}, a multi-pattern is a pattern of the form
$\tau=\tau^0\mn\tau^1\mn\cdots\mn\tau^s$, where 
$\{\tau^0,\tau^1,\dots,\tau^s\}$ is a set of generalized patterns with no hyphens 
and each letter of $\tau^i$ is incomparable with any letter of $\tau^j$ whenever $i\neq j$.

The simplest non-trivial example of a multi-pattern is the multi-pattern $\phi=1\mn1'2'$. To avoid $\phi$ is 
the same as to avoid the patterns $1\mn12$, $1\mn23$, $2\mn12$, $2\mn13$,
and $3\mn12$ simultaneously. To count the number of words in $[k]^n(1\mn1'2')$, 
we choose the leftmost letter of $\sigma \in [k]^n(1\mn1'2')$ in $k$ ways, and observe that all the other letters of $\sigma$ must be in a non-increasing order. Using~\cite{BM1}, for all $n\geq1$, we have
$$|[k]^n(1\mn1'2')|=k\cdot\binom{n+k-2}{n-1}.$$ 

The following theorem is an analogue to~\cite[Theorem~21]{Ki2}.

\begin{theorem}\label{thmpa}
Let $\tau=\tau^0\mn\tau^1$ and $\phi=f_1(\tau^0)\mn f_2(\tau^1)$,
where $f_1$ and $f_2$ are any of the trivial bijections. Then
$\tau\equiv\phi$.
\end{theorem}
\begin{proof}
First, let us prove that the pattern $\tau=\tau^0\mn\tau^1$ is equivalent to the pattern $\phi=\tau^0\mn
f(\tau^1)$, where $f$ is a trivial bijection. Suppose that
$\sigma=\sigma^1\sigma^2\sigma^3\in[k]^n$ avoids $\tau$ and 
$\sigma^1\sigma^2$ has exactly one occurrence of $\tau^0$, namely
$\sigma^2$. Then $\sigma^3$ must avoid $\tau^1$, so $f(\sigma^3)$
avoids $f(\tau^3)$ and $\sigma_f=\sigma^1\sigma^2 f(\sigma^3)$
avoids $\phi$. The converse is also true, if $\sigma_f$ avoids
$\phi$ then $\sigma$ avoids $\tau$. Since any word either avoids
$\tau^0$ or can be factored as above, we have a bijection between
the class of words avoiding $\tau$ and the class of words avoiding $\phi$. 
Thus $\tau\equiv\phi$.

Now, we use the considerations above as well as the properties of
trivial bijections to get
$$\begin{array}{l}
\tau\equiv\tau^0\mn f_2(\tau^1)\equiv R(\tau^0\mn
f_2(\tau^1))\equiv R(f_2(\tau^1))\mn R(\tau^0)\equiv\\
\qquad\qquad\qquad\qquad\qquad\equiv R(f_2(\tau^1))\mn
f_1(R(\tau^0))\equiv R(f_2(\tau^1))\mn R(f_1(\tau^0))\equiv
f_1(\tau^0)\mn f_2(\tau^1).
\end{array}$$
\end{proof}

Using Theorem~\ref{thmpa}, we get the following corollary, which is an analogue 
to~\cite[Corollary 22]{Ki2}.

\begin{corollary}\label{corolperm}
The multi-pattern $\tau^0\mn\tau^1$ is equivalent to the multi-pattern
$\tau^1\mn\tau^0$.
\end{corollary}
\begin{proof}
From Theorem~\ref{thmpa}, using the properties of the trivial bijection $R$, we get
$$\tau^0\mn\tau^1\equiv \tau^0\mn R(\tau^1)\equiv R(R(\tau^1))\mn
R(\tau^0)\equiv \tau^1\mn R(R(\tau^0))\equiv \tau^1\mn\tau^0.$$
\end{proof}

Using induction on $s$, Corollary~\ref{corolperm}, and proceeding in the way proposed 
in~\cite[Theorem 23]{Ki2}, we get

\begin{theorem}\label{thmpb}
Suppose we have multi-patterns $\tau=\tau^0\mn\tau^1\mn\cdots\mn\tau^s$ and
$\phi=\phi^0\mn\phi^1\mn\cdots\mn\phi^s$, where $\tau^1\tau^2\cdots\tau^s$ 
is a permutation of $\phi^1\phi^2\cdots\phi^s$. Then $\tau\equiv\phi$.
\end{theorem}

The last theorem is an analogue to~\cite[Theorem 23]{Ki2}. As a corollary to Theorem~\ref{thmpb}, using Theorem~\ref{thmpa} and the idea of the 
proof of~\cite[Corollary 24]{Ki2}, we get the following corollary which is an analogue to~\cite[Corollary 24]{Ki2}. 

\begin{corollary}
Suppose we have multi-patterns $\tau=\tau^0\mn\tau^1\mn\cdots\mn\tau^s$ and
$\phi=f_0(\tau^0)\mn f_1(\tau^1)\mn\cdots\mn f_s(\tau^s)$, where $f_i$ is 
an arbitrary trivial bijection. Then $\tau\equiv\phi$.
\end{corollary}

The following theorem is a good auxiliary tool for calculating the g.f. for the number of 
words that avoid a given \P. For particular \Ps, it allows to reduce the problem to calculating the g.f. for the number of words that avoid another \P\, which is shorter. We recall that $A_{\tau}^*(x;k)$ is the generating function for the
number of words in $[k]^n$ that quasi-avoid the pattern $\tau$. 

\begin{theorem}\label{thmpc}
Suppose $\tau=\tau^0\mn\phi$, where $\phi$ is an arbitrary \P, and the
letters of $\tau^0$ are incomparable to the letters of $\phi$.
Then for all $k\geq 1$, we have
$$A_\tau(x;k)=A_{\tau^0}(x;k)+A_\phi(x;k)A_{\tau^0}^*(x;k).$$
\end{theorem}
\begin{proof}
Suppose $\sigma=\sigma^1\sigma^2\sigma^3\in [k]^n$ avoids the pattern $\tau$, where $\sigma^1\sigma^2$ quasi-avoids the pattern $\tau^0$, and
$\sigma^2$ is the occurrence of $\tau^0$. Clearly, $\sigma^3$ must avoid $\phi$. To find $A_\tau(x;k)$, we observe that there are two possibilities: 
either $\sigma$ avoids $\tau^0$, or $\sigma$ does not avoid $\tau^0$. In these cases, the 
g.f. for the number of such words is equal to $A_{\tau^0}(x;k)$ and $A_\phi(x;k)A_{\tau^0}^*(x;k)$ 
respectively (the second term came from the factorization above). Thus, the statement is true.
\end{proof}

\begin{corollary}\label{coroldescents}
Let $\tau=\tau^1\mn\tau^2\mn\cdots\mn\tau^s$ be a multi-pattern
such that $\tau^j$ is equal to either $12$ or $21$, for $j=1,2,\ldots,s$. Then
$$A_\tau(x;k)=\frac{1-\left(1+\frac{kx-1}{(1-x)^k}\right)^s}{1-kx}.$$
\end{corollary}
\begin{proof}
According to~\cite{BM2}, $A_{12}(x;k)=A_{21}(x;k)=\frac{1}{(1-x)^k}$. Using Theorem~\ref{thmpc}, 
Proposition~\ref{praa} and induction on $s$, we get the desired result.
\end{proof}

More generally, using Theorem~\ref{thmpc} and Proposition~\ref{praa}, we get the following 
theorem that is the basis for calculating the number of words that avoid a multi-pattern, and therefore is the main result for multi-patterns in this paper.

\begin{theorem}\label{coaa}
Let $\tau=\tau^1\mn\tau^2\mn\cdots\mn\tau^s$ be a multi-pattern.
Then
$$A_\tau(x;k)=\sum_{j=1}^s A_{\tau^j}(x;k)\prod_{i=1}^{j-1} ((kx-1)A_{\tau^i}(x;k)+1).$$
\end{theorem}

\section{The distribution of non-overlapping generalized patterns}

A descent in a word $\sigma\in [k]^n$ is an $i$ such that
$\sigma_i>\sigma_{i+1}$. Two descents $i$ and $j$ {\em overlap} if
$j=i+1$. We define a new statistics, namely the {\em maximum number of
non-overlapping descents,} or $\MND$, in a word. For example,
$\MND(33211)=1$ whereas $\MND(13211143211)=3$. One can find the distribution 
of this new statistic by using Corollary~\ref{coroldescents}. This distribution is 
given in Example~\ref{ex5.2}. However, we prove a more general theorem:

\begin{theorem}\label{dstmnd}
Let $\tau$ be a generalized pattern with no hyphens. Then for all $k\geq 1$,
$$\sum_{n\geq0}\sum_{\sigma\in [k]^n} y^{N_{\tau}(\sigma)}x^n=\frac{A_\tau(x;k)}{1-y((kx-1)A_\tau(x;k)+1)},$$
where $N_{\tau}(\sigma)$ is the maximum number of non-overlapping occurrences of $\tau$ in $\sigma$.
\end{theorem}
\begin{proof}
We fix the natural number $s$ and consider the multi-pattern
$\Phi_s=\tau\mn\tau\mn\cdots\mn\tau$ with $s$ copies of $\tau$. If
a word avoids $\Phi_s$ then it has at most $s-1$ non-overlapping
occurrences of $\tau$. Theorem~\ref{coaa} yields
$$A_{\Phi_s}(x;k)=\sum_{j=1}^s A_{\tau}(x;k)\prod_{i=1}^{j-1} ((kx-1)A_{\tau}(x;k)+1).$$
So, the g.f. for the number of words that has exactly $s$ non-overlapping occurrences of the pattern
$\tau$ is given by
$$A_{\Phi_{s+1}}(x;k)-A_{\Phi_s}(x;k)=A_\tau(x;k)((kx-1)A_\tau(x;k)+1)^s.$$
Hence,
$$\sum_{n\geq0}\sum_{\sigma\in [k]^n}y^{N_{\tau}(\sigma)}x^n=\sum_{s\geq0}
A_\tau(x;k)((kx-1)A_\tau(x;k)+1)^s=\frac{A_\tau(x;k)}{1-y((kx-1)A_\tau(x;k)+1)}.$$
\end{proof}

All of the following examples are corollaries to Theorem~\ref{dstmnd}.

\begin{example}\label{ex5.2}
If we consider descents {\rm(}the pattern $12${\rm)} then $A_{12}(x;k)=\frac{1}{(1-x)^k}$ {\rm(}see~\cite{BM2}{\rm)}, hence the distribution of $\MND$ is given by the formula: 
$$\sum_{n\geq0}\sum_{\sigma\in [k]^n}y^{N_{12}(\sigma)}x^n
=\frac{1}{(1-x)^k+y(1-kx-(1-x)^k)}.$$
\end{example}

\begin{example} The distribution of the maximum number of non-overlapping occurrences of the pattern $122$ is given by the formula:
$$\sum_{n\geq0}\sum_{\sigma\in [k]^n}y^{N_{122}(\sigma)}x^n
=\frac{x}{(1-x^2)^k+x-1+y(1-kx^2-(1-x^2)^k)},$$
since according to~\cite[Theorem 3.10]{BM3}, $A_{122}(x;k)=\frac{x}{(1-x^2)^k-(1-x)}$.
\end{example}

\begin{example} If we consider the pattern $212$ then $A_{212}(x;k)=\left(1-x\displaystyle\sum_{j=0}^{k-1}\frac{1}{1+jx^2}\right)^{-1}$ {\rm(}see~\cite[Theorem 3.12]{BM3}{\rm)}, hence the distribution of the maximum number of non-overlapping occurrences of the pattern $212$ is given by the formula:
$$\sum_{n\geq0}\sum_{\sigma\in [k]^n}y^{N_{212}(\sigma)}x^n
=\frac{1}{1-x\displaystyle\sum_{j=0}^{k-1}\frac{1}{1+jx^2}+xy\left(\displaystyle\sum_{j=0}^{k-1}\frac{1}{1+jx^2}-k\right)}.$$
\end{example}

\begin{example} Using~\cite[Theorem 3.13]{BM3}, the distribution of the maximum number of non-overlapping occurrences of the pattern $123$ is given by the formula:
$$\sum_{n\geq0}\sum_{\sigma\in [k]^n}y^{N_{123}(\sigma)}x^n
= \frac{1}{\displaystyle\sum_{j=0}^{k}a_j{k \choose j}x^j + y\left(1-kx-\displaystyle\sum_{j=0}^{k}a_j{k \choose j}x^j \right)},$$
where $a_{3m}=1$, $a_{3m+1}=-1$, and $a_{3m+2}=0$, for all $m\geq 0$. 
\end{example}



\begin{thebibliography}{WWW}
\bibitem[BS]{BS} E. Babson, E. Steingr\'{\i}msson: Generalized permutation
patterns and a classification of the Mahonian statistics, S\'eminaire
Lotharingien de Combinatoire, B44b:18pp, (2000).

\bibitem[Bu]{Burstein}
A.~Burstein, Enumeration of words with forbidden patterns, Ph.D.
thesis, University of Pennsylvania, 1998.

\bibitem[BM1]{BM1}
A.~Burstein and T.~Mansour, Words restricted by patterns with at
most 2 distinct letters, \emph{Electronic J. of Combinatorics}, to
appear (2002).

\bibitem[BM2]{BM2}
A.~Burstein and T.~Mansour, Words restricted by $3$-letter
generalized multipermutation patterns, preprint CO/0112281.

\bibitem[BM3]{BM3}
A.~Burstein and T.~Mansour, Counting occurrences of some subword
patterns, preprint CO/0204320.

\bibitem[C]{Claesson}
A. Claesson: Generalised Pattern Avoidance, European Journal of Combinatorics {\bf 22} (2001), 961-971.

\bibitem[CM]{CM}
A.~Claesson and T.~Mansour, Enumerating Permutations Avoiding a Pair of Babson-Steingr\'{\i}msson Patterns, preprint CO/0107044.

\bibitem[Ki1]{Ki1}
S.~Kitaev, Multi-avoidance of generalised patterns, Discr. Math., to appear (2002).

\bibitem[Ki2]{Ki2}
S.~Kitaev, Partially ordered generalized patterns, Discr. Math., to appear (2002).

\bibitem[Kn]{Knuth}
\bibitem[Knuth]{Knuth} D. E. Knuth: {\em The Art of Computer Programming}, 2nd ed. Addison Wesley, Reading, MA, (1973).

\bibitem[SS]{SS}
R. Simion, F. Schmidt: Restricted permutations, European
J. Combin. {\bf 6}, no. 4 (1985), 383--406.
\end{thebibliography}
\end{document}